\documentclass[reqno]{amsart}
\usepackage{amssymb}
\usepackage{amsmath}

\makeatletter
\@addtoreset{equation}{section}
\makeatother

\renewcommand\thefigure{\thesection.\@arabic\c@figure}
\renewcommand\thetable{\thesection.\@arabic\c@table}
 
\newtheorem{theorem}{Theorem}[section]
\newtheorem{lemma}[theorem]{Lemma}
\newtheorem{proposition}[theorem]{Proposition}

\newcommand{\mc}[1]{{\mathcal #1}}

\newcommand{\bb}[1]{{\mathbb #1}}

\renewcommand{\epsilon}{\varepsilon}
\renewcommand{\tilde}{\widetilde}

\newcommand{\<}{\langle}
\renewcommand{\>}{\rangle}
\newcommand{\Bl}{\Big\langle}
\newcommand{\Br}{\Big\rangle}

\begin{document}

\author{L.\ Bertini, A.\ De Sole, D.\ Gabrielli, G.\ Jona-Lasinio,
C.\ Landim} 

\address{
Lorenzo Bertini 
\newline
Dipartimento di Matematica, Universit\`a di Roma ``La Sapienza'' 
\newline
Piazzale Aldo Moro 2, 00185 Roma, Italy
\newline
e-mail: \rm \texttt{bertini@mat.uniroma1.it}
}

\address{
Alberto De Sole
\newline
Department of Mathematics, Harvard University, 
\newline
Cambridge, MA 02138, USA
\newline
e-mail:  \rm \texttt{desole@math.harvard.edu}}

\address{
Davide Gabrielli
\newline
Dipartimento di Matematica, Universit\`a dell'Aquila,
\newline
Coppito, 67100 L'Aquila, Italy
\newline
e-mail:  \rm \texttt{gabriell@univaq.it}}

\address{
Giovanni Jona-Lasinio
\newline
Dipartimento di Fisica and INFN, Universit\`a di Roma ``La Sapienza", 
\newline
Piazzale A. Moro 2, 00185 Roma, Italy 
\newline 
e-mail:  \rm \texttt{gianni.jona@roma1.infn.it}}

\address{
Claudio Landim
\newline
IMPA, Estrada Dona Castorina 110, CEP 22460 Rio de Janeiro, Brasil 
\newline
and CNRS UMR 6085,
Universit\'e de Rouen, 76128 Mont Saint Aignan, France.
\newline
e-mail:  \rm \texttt{landim@impa.br}}

\title[Dynamical large deviations in stochastic lattice gases]
{Large deviation approach to non equilibrium processes in stochastic 
lattice gases} 

\begin{abstract}
We present a review of recent work on the statistical mechanics of non
equilibrium processes based on the analysis of large deviations 
properties of microscopic systems. Stochastic lattice gases are 
non trivial models of such phenomena and can be studied rigorously    
providing a source of challenging mathematical problems. 
In this way, some principles of wide validity have been obtained
leading to interesting physical consequences. 
\end{abstract}

                                           

\maketitle

\section{A Physicist motivation}
\label{sec0}     

In equilibrium statistical mechanics there is a well defined
relationship, established by Boltzmann, between the probability of a
state and its entropy. This fact was exploited by Einstein to study
thermodynamic fluctuations.  So far it does not exist a theory of
irreversible processes of the same generality as equilibrium
statistical mechanics and presumably it cannot exist. While in
equilibrium the Gibbs distribution provides all the information and no
equation of motion has to be solved, the dynamics plays the major role
in non equilibrium.

When we are out of equilibrium, for example in a stationary state of a
system in contact with two reservoirs, even if the system is in a
local equilibrium state so that it is possible to define the local
thermodynamic variables e.g.\ density or magnetization, it is not
completely clear how to define the thermodynamic potentials like the
entropy or the free energy.  One possibility, adopting the
Boltzmann-Einstein point of view, is to use fluctuation theory to
define their non equilibrium analogs. In fact, in this way extensive
functionals can be obtained although not necessarily simply additive
due to the presence of long range correlations which seem to be a
rather generic feature of non equilibrium systems.

Let us recall the Boltzmann-Einstein theory of equilibrium
thermodynamic fluctuations. The main principle is that the probability
of a fluctuation in a macroscopic region of fixed volume $V$ is
\begin{equation}
\label{1}
P \propto \exp\{V\Delta S / k\}
\end{equation}
where $\Delta S$ is the variation of the specific entropy calculated
along a reversible transformation creating the fluctuation and $k$ is
the Boltzmann constant.  Eq.\ (\ref{1}) was derived by Einstein simply
by inverting the Boltzmann relationship between entropy and
probability.  He considered (\ref{1}) as a phenomenological definition
of the probability of a state.  Einstein theory refers to fluctuations
for equilibrium states, that is for systems isolated or in contact
with reservoirs characterized by the same chemical potentials.  When
in contact with reservoirs $\Delta S$ is the variation of the total
entropy (system + reservoirs) which for fluctuations of constant
volume and temperature is equal to $-{\Delta \mc F}/{T}$, that is minus
the variation of the free energy of the system divided by the
temperature.

We consider a stationary non-equilibrium state (SNS), namely, due to
external fields and/or different chemical potentials at the
boundaries, there is a flow of physical quantities, such as heat,
electric charge, chemical substances, across the system.  To start
with, it is not always clear that a closed macroscopic dynamical
description is possible.  If the system can be described by a
hydrodynamic equation, a fact which can be rigorously established in
stochastic lattice gases, a reasonable goal is to find an explicit
connection between the thermodynamic potentials and the dynamical
macroscopic properties like transport coefficients.  
The study of large fluctuations provides such a connection. 

Besides the definition of thermodynamic potentials, in a dynamical
setting a typical question one may ask is the following: what is the
most probable trajectory followed by the system in the spontaneous
emergence of a fluctuation or in its relaxation to an equilibrium or a
stationary state?  To answer this question one first derives a
generalization of the Boltzmann-Einstein formula from which the most
probable trajectory can be calculated by solving a variational
principle.  
For equilibrium states and small fluctuations an answer to this type of
questions was given by Onsager and Machlup in 1953 \cite{on}.  The
Onsager-Machlup theory gives the following result under the assumption
of time reversibility of the microscopic dynamics:
the most probable creation and relaxation trajectories of a
fluctuation are one the time reversal of the other.

We discuss this issue in the context of stochastic lattice gases in a
box of linear size $N$ with birth and death process at the boundary
modeling the reservoirs.  
We consider the case when there is only one thermodynamic variable, 
the local density denoted by  $\rho$. 
Its macroscopic evolution is given by the continuity equation 
\begin{equation}
\label{H1}
\partial_t \rho  \;=\;
\nabla \cdot \Big[ D(\rho)\nabla \rho - \chi(\rho) E \Big]
\;=\; - \nabla\cdot J(\rho)
\end{equation}
where $D(\rho)$ is the diffusion matrix, $\chi(\rho)$ the mobility and
$E$ the external field. Here $J(\rho)$ is the macroscopic
instantaneous current associated to the density profile $\rho$.
Finally the interaction with the reservoirs
appears as boundary conditions to be imposed on solutions of
(\ref{H1}). We shall denote by $u$ the macroscopic space coordinate
and by $\bar\rho=\bar\rho(u)$ the unique stationary
solution of (\ref{H1}), i.e.\ $\bar\rho$ is the typical density profile
for the SNS.

This equation derives from the underlying stochastic dynamics through
an appropriate scaling limit in which the microscopic time and space
coordinates are rescaled diffusively. The hydrodynamic equation
(\ref{H1}) thus represents the law of large numbers for the empirical
density of the stochastic lattice gas. The convergence has to be
understood in probability with respect to the law of the stochastic
lattice gas. Finally, the initial condition for \eqref{H1} depends on
the initial distribution of particles. Of course many microscopic
configurations give rise to the same initial condition $\rho_0(u)$.

Let us denote by $\nu^N$ the invariant measure of the stochastic
lattice gas. The free energy $\mc F(\rho)$, defined as a functional of the
density profile $\rho=\rho(u)$, gives the asymptotic probability
of fluctuations of the empirical measure $\pi^N$ under 
the invariant measure $\nu^N$. More precisely 
\begin{equation}
\label{BE2}
\nu^N \big( \pi^N \approx \rho \big) \sim  
\exp\big\{ - N^d \mc F(\rho) \big\}
\end{equation}
where $d$ is the dimensionality of the system, $\pi^N \approx \rho$ 
means closeness in the weak topology and $\sim$ denotes
logarithmic equivalence as $N\to\infty$. 
In the above formula we omitted the dependence on the
temperature since it does not play any role in our analysis; we also 
normalized $\mc F$ so that $\mc F(\bar\rho)=0$. 

In the same way, the behavior of space time fluctuations can be
described as follows. Let us denote by $\bb P_{\nu^N}$ the stationary
process of the stochastic lattice gas, i.e.\ the initial distribution
is given by the invariant measure $\nu^N$.
The probability that the
evolution of the random variable $\pi^N_t$ deviates from the solution
of the hydrodynamic equation and is close to some trajectory
${\hat{\rho}}_t$ is exponentially small  and of the form
\begin{equation}
\label{LD}
\bb P_{\nu^N} \big( \pi^N_t \approx  \hat{\rho}_t, \: 
t\in [t_1, t_2]\big) \sim
\exp\big\{-N^d 
\big[ \mc F(\hat \rho_{t_1}) + {I}_{[t_1,t_2]}(\hat{\rho}) 
\big] \big\}
\end{equation}
where $I(\hat{\rho})$ is a functional which vanishes if
${\hat{\rho}}_t$ is a solution of (\ref{H1}) and $\mc F(\hat
\rho_{t_1})$ is the free energy cost to produce the initial density
profile ${\hat{\rho}}_{t_1}$.  Therefore $I(\hat{\rho})$
represents the extra cost necessary to follow the trajectory
${\hat{\rho}}_t$ in the time interval $[t_1,t_2]$.  

To determine the most probable trajectory followed by the system in
the spontaneous creation of a fluctuation, we consider the following
physical situation.  The system is macroscopically in the stationary
state $\bar\rho$ at $t=-\infty$ but at $t=0$ we find it in the state
$\rho$.  According to (\ref{LD}) the most probable trajectory is the
one that minimizes $I$ among all trajectories $\hat\rho_t$
connecting $\bar\rho$ to ${\rho}$ in the time interval
$[-\infty,0]$, that is the optimal path for the variational problem 
\begin{equation}
\label{l1}
V(\rho)= \inf_{\hat \rho}  I_{[-\infty,0]}(\hat \rho)
\end{equation}
The functional $V(\rho)$, called the quasi-potential, measures the
probability of the fluctuation $\rho$.  Moreover, the optimal
trajectory for \eqref{l1} determines the path followed by the system
in the creation of the fluctuation $\rho$.  As shown in
\cite{bdgjl1,bdgjl2,bg} this minimization problem gives the non
equilibrium free energy, i.e.\ $V=\mc F$.
As we discuss here, by analyzing this variational problem 
for SNS, the Onsager-Machlup relationship has to be modified in the
following way: the spontaneous emergence of a macroscopic
fluctuation takes place most likely following a trajectory which can
be characterized in terms of the time reversed process.

\bigskip
Beside the density, a very important observable is the current flux. 
This quantity
gives informations that cannot be recovered from the density because
from a density trajectory we can determine the current trajectory only
up to a divergence free vector field.  We emphasize that this is due to the
loss of information in the passage from the microscopic level to the
macroscopic one.

To discuss the current fluctuations in the context of stochastic
lattice gases, we introduce the empirical current $w^N$ which measures
the local net flow of particles.  As for the empirical density, it is
possible to prove a dynamical large deviations principle for the
empirical current which is informally stated as follow.
Given a vector field $j:[0,T]\times \Lambda \to \bb R^d$, we have
\begin{equation}
\label{f1a}
\bb P_{\eta^N} \big( w^N  \approx j (t,u) \big) 
\sim \exp\big\{ - N^d \, \mc I_{[0,T]}(j)\big\}
\end{equation}
where $\bb P_{\eta^N}$ is the law of the stochastic lattice gas with
initial condition given by $\eta^N=\{\eta^N_x\}$, which represents 
the number of  particles in each site, and the rate functional is 
\begin{equation}
\label{Ica}
\mc I_{[0,T]}(j)\;=\; \frac 12 \int_0^T \!dt \,
\big\langle [ j - J(\rho)  ], \chi(\rho)^{-1}
[ j - J(\rho) ] \big\rangle
\end{equation}
in which we recall that 
\begin{equation*}
J(\rho) = - D(\rho) \nabla \rho  + \chi(\rho) E\; .
\end{equation*}
Moreover, $\rho=\rho(t,u)$ is obtained by solving the continuity
equation $\partial_t\rho +\nabla\cdot j =0$ with the initial condition
$\rho(0)=\rho_0$ associated to $\eta^N$. The rate functional vanishes if
$j=J(\rho)$, so that $\rho$ solves \eqref{H1}. This is the law of
large numbers for the observable $w^N$. Note that equation
\eqref{Ica} can be interpreted, in analogy to the classical Ohm's law,
as the total energy dissipated in the time interval $[0,T]$ by the
extra current $j - J(\rho)$.

Among the many problems we can discuss within this theory, we study
the fluctuations of the time average of the empirical current over a
large time interval.  
We show that the probability of observing a time-averaged 
fluctuation $J$ can be described by a functional $\Phi(J)$
which we characterize in terms of a variational
problem for the functional $\mc I_{[0,T]}$
\begin{equation}
\label{limTa}
\Phi (J)
 = \lim_{T\to\infty} \; \inf_{j} 
\frac 1T \; \mc I_{[0,T]} (j)\;,
\end{equation}
where the infimum is carried over all paths $j=j(t,u)$ having time
average $J$.  We finally analyze the variational problem \eqref{limTa}
for some stochastic lattice gas models and show that different
scenarios take place. In particular, for the symmetric exclusion
process with periodic boundary condition the optimal trajectory is
constant in time. On the other hand for the KMP model \cite{kmp}, also
with periodic boundary conditions, this is not the case: we show that
a current path in the form of a traveling wave leads to a higher
probability.

\section{Boundary driven simple exclusion process}
\label{sec1}     

For an integer $N\ge 1$, let $\Lambda_N := \{1 ,\dots , N-1\}$.  The
sites of $\Lambda_N$ are denoted by $x$, $y$, and $z$ while the
macroscopic space variable (points in the interval $[0,1]$) by $u$. 
We introduce the microscopic state space as $\Sigma_N:=
\{0,1\}^{\Lambda_N}$ which is endowed with the discrete topology;
elements of $\Sigma_N$, called configurations, are denoted by $\eta$.
In this way $\eta (x)\in \{0,1\}$ stands for the number of particles
at site $x$ for the configuration $\eta$.

The one dimensional boundary driven simple exclusion process is the
Markov process on the state space $\Sigma_N$ with infinitesimal
generator defined as follows. Given $\alpha,\beta \in (0,1)$ we let
\begin{eqnarray*}
(L_{N} f)(\eta) &:=&  \frac {N^2}2 \sum_{x=1}^{N-2}  
\big[ f(\sigma^{x,x+1}\eta) - f(\eta) \big] \\
&+ &  \frac {N^2}2\,  [\alpha \{1-\eta(1)\} + (1-\alpha) \eta(1)] 
\big[ f(\sigma^{1} \eta) - f(\eta) \big]   \\
&+ &  \frac {N^2}2\,  [\beta \{1-\eta(N-1)\} + (1-\beta) \eta(N-1)] 
\big[ f(\sigma^{N-1} \eta) - f(\eta) \big]   
\end{eqnarray*}
for every function $f: \Sigma_N\to\bb R$.  In this formula
$\sigma^{x,y}\eta$ is the configuration obtained from $\eta$ by
exchanging the occupation variables $\eta (x)$ and $\eta (y)$:
\begin{equation*}
(\sigma^{x,y} \eta) (z) \; :=\; 
\left\{\begin{array}{ll}
                \eta (y) & {\rm if \ } z=x\\
                \eta (x) & {\rm if \ } z=y\\
                \eta (z)   & {\rm if \ } z\neq x,y
                \end{array}
\right.  
\end{equation*}
and $\sigma^{x}\eta$ is the configuration obtained from $\eta$ by
flipping the configuration at $x$: 
\begin{equation*}
(\sigma^x \eta)\,(z) \; :=\;
\eta (z) [1-\delta_{x,z}] \;+\; \delta_{x,z} [1-\eta (z)]\; ,
\end{equation*}
where
$\delta_{x,y}$ is the Kronecker delta.  
The parameters $\alpha,\beta$, which affect the birth and death rates
at the two boundaries, represent the densities of the
reservoirs. Without loss of generality, we assume $\alpha \le \beta$.
Notice finally that $L_N$ has been speeded up by $N^2$; this
corresponds to the diffusive scaling.

The Markov process $\{\eta_t : t\ge 0\}$ associated to the generator
$L_N$ is irreducible. It has therefore a unique invariant measure, 
denoted by $\nu^N_{\alpha, \beta}$. 
The process is reversible if and only if $\alpha = \beta$, in which
case $\nu^N_{\alpha, \alpha}$ is the
Bernoulli product measure with density $\alpha$
\begin{equation*}
\nu^N_{\alpha, \alpha} \{\eta : \eta(x) = 1\} \;=\; \alpha
\end{equation*}
for $1\le x\le N-1$.

If $\alpha \not = \beta$ the process is  not reversible and the
measure $\nu^N_{\alpha, \beta}$ carries long range correlations. 
Since $E_{\nu^N_{\alpha, \beta}} [L_N \eta(x)] =
0$, it is not difficult to show that $\rho^N(x) = E_{\nu^N_{\alpha,
    \beta}} [\eta(x)]$ is the solution of the linear equation
\begin{equation}
\label{f01}
\left\{
\begin{array}{l}
\Delta_N \rho^N(x) = 0 \; , \quad 1\le x\le N-1\;, \\
\rho^N(0) = \alpha\; , \quad \rho^N(N) = \beta\;,
\end{array}
\right.
\end{equation}
where $\Delta_N$ stands for the discrete Laplacian.  
Hence 
\begin{equation}
\label{f01bis}
\rho^N(x) = \alpha  + \frac xN \, (\beta -\alpha)
\end{equation}

It is also possible to obtain a closed expression for the correlations
$$
E_{\nu^N_{\alpha, \beta}} [\eta(x); \eta(y)] 
= E_{\nu^N_{\alpha,\beta}} [\eta(x) \eta(y)] - 
E_{\nu^N_{\alpha, \beta}} [\eta(x)] 
E_{\nu^N_{\alpha, \beta}} [\eta(y)]
$$ 
As shown in \cite{DFIP,S}, for $1\le x< y\le N-1$ we have
\begin{equation}
\label{f03}
E_{\nu^N_{\alpha, \beta}} [\eta(x); \eta(y)]
= - \frac { (\beta -\alpha)^2}{N-1} \, \frac xN  \Big( 1 - \frac yN \Big)
\end{equation}
Computing $L_N \eta(x) \eta(y)$, we obtain that the correlations solve 
a discrete differential equation. One can then check that 
\eqref{f03} is the solution.

Note that, if we take $x,y$ at distance $O(N)$ from the boundary, then
the covariance between $\eta(x)$ and $\eta(y)$ is of order $O(1/N)$.
Moreover the random variables $\eta(x)$ and $\eta(y)$ are negatively
correlated. This is the same qualitative behavior as the one in the 
canonical Gibbs measure given by the uniform measure on 
$\Sigma_{N,k} =\{\eta\in\Sigma_N \, :\, \sum_{x=1}^{N-1}\eta(x)=k\}$.

\section{Stationary large deviations of the empirical density}
\label{sec2}

Denote by $\mc M_+$ the space of positive measures on $[0,1]$ with
total mass bounded by $1$. We consider $\mc M_+$ endowed with the weak
topology. For a configuration $\eta$ in $\Sigma_N$, let $\pi^N$ be the
measure obtained by assigning mass $N^{-1}$ to each particle and
rescaling space by $N^{-1}$
\begin{equation*}
\pi^N (\eta) \;:=\; \frac 1N \sum_{x=1}^{N-1} \eta (x) \, \delta_{x/N}\;,
\end{equation*}
where $\delta_u$ stands for the Dirac measure concentrated on $u$.
Denote by $\<\pi^N, H\>$ the integral of a continuous function
$H:[0,1]\to \bb R$ with respect to $\pi^N$
\begin{equation*}
\<\pi^N, H\> \;=\; \frac 1N \sum_{x=1}^{N-1} H(x/N) \eta(x)\;.
\end{equation*}
We use the same notation for the inner product in
$L_2([0,1],du)$. Analogously we denote the space integral of a function $f$
by $\<f\>=\int_0^1\!du\, f(u)$.

The law of large numbers for the empirical density under the
stationary state $\nu^N_{\alpha, \beta}$ is proven in 
\cite{DFIP,els1,els2}.

\begin{theorem}
\label{s1}
For every continuous function $H:[0,1]\to \bb R$ and every $\delta>0$, 
\begin{equation*}
\lim_{N\to\infty} \nu^N_{\alpha, \beta} \Big\{ \, \Big\vert \<\pi^N,
H\> - \<\bar\rho , H \> \Big\vert > \delta \Big\} \;=\;
0\;, 
\end{equation*}
where 
\begin{equation}
\label{f02}
\bar \rho(u) = \alpha (1-u) +\beta u\,.
\end{equation}
\end{theorem}

We remark that $\bar\rho$ is the solution of the elliptic linear equation 
$$
\left\{
\begin{array}{l}
\Delta \rho = 0 \; ,  \\
\rho(0) = \alpha\; , \quad \rho(1) = \beta\;.
\end{array}
\right.
$$
which is the continuous analog of \eqref{f01}. 
Here and further $\Delta$ stands for the Laplacian.

Once a law of large numbers has been established, it is natural to
consider the deviations around the typical value $\bar \rho$. 
From the explicit expression of the microscopic correlations
\eqref{f03} it is possible to prove a central limit theorem for the
empirical density under the stationary measure $\nu^N_{\alpha,\beta}$. 
We refer to \cite{S} for  a more detailed discussion and to 
\cite{GKS} for the mathematical details.

Fix a profile $\gamma: [0,1] \to [0,1]$ different from $\bar\rho$ and
a neighborhood $V_\epsilon(\gamma)$ of radius $\epsilon>0$ around the
measure $\gamma(u) du$ in $\mc M_+$.  The mathematical formulation of
the Boltzmann-Einstein formula \eqref{1}  consists in
determining the exponential rate of decay, as $N\uparrow\infty$, of
\begin{equation*}
\nu^N_{\alpha, \beta} \big\{ \pi^N \in V_\epsilon(\gamma) \big\}\;.
\end{equation*}
Derrida, Lebowitz and Speer \cite{dls1, dls2} derived, by explicit
computations, the large deviations principle for the empirical density
under the stationary state $\nu^N_{\alpha, \beta}$. 
This result has been obtained by a dynamical/variational approach in 
\cite{bdgjl2}, a rigorous proof is given in \cite{bdgjl3}. The precise
statement is the following.

\begin{theorem}
\label{s02}
For each profile $\gamma: [0,1]\to [0,1]$,
\begin{eqnarray*}
\!\!\!\!\!\!\!\!\!\!\!\!\!\!\!\!\!\! &&
\limsup_{\epsilon\to 0} \limsup_{N\to\infty} \frac 1N \log
\nu^N_{\alpha, \beta} \big\{ \pi^N \in V_\epsilon(\gamma) \big\}
\;\le\; - \mc F (\gamma) \; , \\
\!\!\!\!\!\!\!\!\!\!\!\!\!\!\!\!\!\! && \quad
\liminf_{\epsilon\to 0} \liminf_{N\to\infty} \frac 1N \log
\nu^N_{\alpha, \beta} \big\{ \pi^N \in V_\epsilon(\gamma) \big\}
\;\ge\; - \mc F (\gamma) \; ,
\end{eqnarray*}
where
\begin{equation}
\label{f11}
\mc F (\gamma) \;=\; \int_{0}^1 \! du\, \Big\{ \gamma (u) \log 
\frac {\gamma (u)} {F(u)} \;+\; [1 - \gamma (u)] 
\log \frac {1- \gamma (u)}{1- F(u)} + 
\log \frac {F'(u)}{\beta - \alpha} \Big\}
\end{equation}
and $F\in C^1([0,1])$ is the unique increasing solution of 
the non linear boundary value problem 
\begin{equation}\label{Deq}
\left\{
\begin{array}{l}
{\displaystyle F'' = \big( \gamma - F \big) 
\frac{\big( F' \big)^2}{F(1-F)} }\; , \\
{\displaystyle F(0) = \alpha\;, \quad F(1) = \beta\;. }
\end{array}
\right.
\end{equation}
\end{theorem}

It is interesting to compare the large deviation properties of the
stationary state $\nu^N_{\alpha, \beta}$ with the one of 
$\mu^N_{\alpha, \beta}$, the product measure on $\Sigma_N$
which has the same marginals as $\nu^N_{\alpha, \beta}$, i.e.\
\begin{equation*}
\mu^N_{\alpha, \beta} \{\eta : \eta(x) = 1\} \;=\; \rho^N(x)\;,
\end{equation*}
where $\rho^N$ is given by \eqref{f01bis}. It is not difficult to
show that in this case
\begin{eqnarray*}
\!\!\!\!\!\!\!\!\!\!\!\!\!\!\!\!\!\! &&
\limsup_{\epsilon\to 0} \limsup_{N\to\infty} \frac 1N \log
\mu^N_{\alpha, \beta} \big\{ \pi^N \in V_\epsilon(\gamma) \big\}
\;\le\; - \mc F_0(\gamma) \; , \\
\!\!\!\!\!\!\!\!\!\!\!\!\!\!\!\!\!\! && \quad
\liminf_{\epsilon\to 0} \liminf_{N\to\infty} \frac 1N \log
\mu^N_{\alpha, \beta} \big\{ \pi^N \in V_\epsilon(\gamma) \big\}
\;\ge\; - \mc F_0(\gamma) \; ,
\end{eqnarray*}
where
\begin{equation*}
\mc F_0(\gamma) \;=\; \int_{0}^1 \! du \,
\Big\{ \gamma (u) \log \frac {\gamma (u)}{\bar\rho (u)} 
\;+\; [1 - \gamma (u)] \log \frac {1- \gamma (u)} {1-\bar\rho (u)} 
\Big\}
\end{equation*}
and $\bar\rho$ is given in \eqref{f02}.  Notice that the functional
$\mc F_0$ is local while $\mc F$ is not.  Moreover, it is not
difficult to show \cite{bdgjl3,dls2} that $\mc F_0\le \mc F$.
Therefore, fluctuations have less probability for the stationary state
$\nu^N_{\alpha,\beta}$ than for the product measure $\mu^N_{\alpha,
\beta}$.  This bound reflects at the large deviations level the
negative correlations observed in \eqref{f03}.

\section{Hydrodynamics and dynamical large deviations of the density}
\label{sec3}

We discuss the asymptotic behavior, as $N\to\infty$, 
of the evolution of the empirical density. 
Denote by $\{\eta^N_t : t\ge 0\}$ the Markov process with generator
$L_N$ and $\pi^N_t = \pi^N(\eta^N_t)$.  Fix a profile $\gamma:
[0,1] \to [0,1]$ and assume that $\pi^N_0$ converges to $\gamma(u) du$ as
$N\uparrow\infty$. Observing the time evolution of the process, we
expect $\pi^N_t$ to relax to the stationary profile
$\bar\rho(u) du$ according to some trajectory $\rho_t(u) du$. 
This result is usually referred to as  the hydrodynamic limit.
For the boundary driven simple exclusion process it is stated in
Theorem~\ref{s03} below \cite{els1, els2}.

Fix $T>0$ and denote, respectively, by $D([0,T], \mc M_+)$, 
$D([0,T], \Sigma_N)$ the space of $\mc M_+$-valued, $\Sigma_N$-valued
cadlag functions endowed with the Skorohod topology.  For a
configuration $\eta^N$ in $\Sigma_N$, denote by $\bb P_{\eta^N}$ the
probability on the path space $D([0,T], \Sigma_N)$ induced by the
initial state $\eta^N$ and the Markov dynamics associated to the
generator $L_N$.

\begin{theorem}
\label{s03}
Fix a profile $\gamma: [0,1]\to [0,1]$  and a sequence of configurations
$\eta^N$ such that $\pi^N(\eta^N)$ converges to $\gamma (u) du$, as
$N\uparrow\infty$. Then, for each $t\ge 0$, $\pi^N_t$ converges in 
$\bb P_{\eta^N}$-probability to $\rho_t(u) du$  as $N\uparrow\infty$. 
Here $\rho_t(u)$ is the solution of the heat equation 
\begin{equation}
\label{f04}
\left\{
\begin{array}{l}
{\displaystyle \partial_t \rho_t(u) = (1/2) \Delta \rho_t(u) \; ,} \\
{\displaystyle \rho_0(u) = \gamma(u) \;, \quad u\in (0,1)} \\
{\displaystyle \rho_t(0) = \alpha\;, 
\quad  \rho_t(1) = \beta\;. }
\end{array}
\right.
\end{equation}
In other words, for each $\delta, T >0$ and each 
continuous function $H:[0,1]\to \bb R$ we have
$$
\lim_{N\to\infty} 
\bb P_{\eta^N} \Big( \sup_{t\in[0,T]} 
\big| \langle \pi^N_t, H \rangle - \langle \rho_t , H \rangle \big| >
\delta \Big) =0
$$
\end{theorem}

Equation \eqref{f04} describes the relaxation path from $\gamma$ to
$\bar \rho$ since $\rho_t$ converges to the stationary path
$\bar\rho$ as $t\uparrow\infty$. To examine the fluctuations paths, we
need first to describe the large deviations of the trajectories in a
fixed time interval. This result requires some notation.

Fix a profile $\gamma$ bounded away from $0$ and $1$:  
for some $\delta>0$ we have $\delta \le \gamma \le 1-\delta$ $du$-a.e.
Denote by $C_\gamma$
the following subset of $D([0,T], \mc M_+)$. A trajectory $\pi_t$,
$t\in[0,T]$ is in $C_\gamma$ if it is continuous and, for any
$t\in [0,T]$, we have $\pi_t (du) = \lambda_t(u) du$ for some density
$\lambda_t (u) \in [0,1]$ which satisfies the boundary conditions
$\lambda_0 = \gamma$, $\lambda_t(0) = \alpha$,
$\lambda_t(1) = \beta$. The latter are to be understood in the
sense that, for each $t\in[0,T]$,
$$
\lim_{\delta\downarrow 0} \frac 1\delta \, 
\int_0^\delta \!du \:  \lambda_t(u) = \alpha\,,
\qquad\qquad
\lim_{\delta\downarrow 0} \frac 1\delta \, 
\int_{1-\delta}^1 \!du \: \lambda_t(u) = \beta\,.
$$

We define a functional $I_{[0,T]} (\cdot |
\gamma)$ on $D([0,T], \mc M_+)$ by setting 
$I_{[0,T]} (\pi | \gamma) = + \infty$ if $\pi\not\in C_\gamma$ 
and by a variational expression for $\pi \in C_\gamma$. Referring to 
\cite[Eq.\ (2.4)--(2.5)]{bdgjl3} for the precise definition, here we
note that if $\pi_t (du) = \lambda_t(u) du$ for some smooth density
$\lambda$ we have 
\begin{equation}
\label{f05}
I_{[0,T]} (\pi | \gamma) \; =\; \frac 12 \int_0^T \!dt 
\int_{0}^1 \!du\, \chi(\lambda_t(u)) [\nabla H_t(u)]^2\; .
\end{equation}
Here, $\chi(a) = a(1-a)$ is the mobility in the symmetric simple
exclusion process and $H_t$ is the unique solution of
\begin{equation}
\label{f07}
\partial_t \lambda_t(u) = (1/2) \Delta \lambda_t(u) -  
\nabla \big[ \chi(\lambda_t(u)) \nabla H_t(u) \big]\, ,
\qquad u\in (0,1) \;
\end{equation}
with the boundary conditions $H_t(0)=H_t(1)=0$ for any $t\in [0,T]$.
As before $\Delta$ stands for $\frac {d^2}{du^2}$ and $\nabla$ 
stands for $\frac {d}{du}$.
Hence, to compute $I_{[0,T]} (\pi | \gamma)$, we first solve
equation \eqref{f07} in $H$ and then plug it in \eqref{f05}.

The following theorem states the dynamical large deviation principle
for the boundary driven simple exclusion process. It is proven in
\cite{bdgjl3} by developing  the techniques introduced in \cite{dv,kov}.

\begin{theorem}
\label{s04}
Fix $T>0$ and a profile $\gamma$ bounded away from $0$ and $1$.
Consider a sequence $\eta^N$ of configurations associated to $\gamma$
in the sense that $\pi^N(\eta^N)$ converges to $\gamma(u) du$ as
$N\uparrow\infty$. Fix $\pi$ in $D([0,T], \mc M_+ )$ and a
neighborhood $V_\epsilon(\pi)$ of $\pi$ of radius $\epsilon$.
Then
\begin{eqnarray*}
\!\!\!\!\!\!\!\!\!\!\!\!\!\!\!\!\!\! &&
\limsup_{\epsilon\to 0} \limsup_{N\to\infty} \frac 1N \log
\bb P_{\eta^N} \big\{ \pi^N \in V_\epsilon(\pi) \big\}
\;\le\; - I_{[0,T]} (\pi | \gamma) \; , \\
\!\!\!\!\!\!\!\!\!\!\!\!\!\!\!\!\!\! && \quad
\liminf_{\epsilon\to 0} \liminf_{N\to\infty} \frac 1N \log
\bb P_{\eta^N} \big\{ \pi^N \in V_\epsilon(\pi) \big\}
\;\ge\; - I_{[0,T]} (\pi | \gamma)  \; .
\end{eqnarray*}
\end{theorem}

We may now formulate the following exit problem. Fix a profile
$\gamma$ and a path $\pi$ such that $\pi_0 = \bar\rho \, du$, $\pi_T =
\gamma \, du$.  The functional $I_{[0,T]} (\pi | \bar\rho)$ measures the
cost of observing the path $\pi$. Therefore,
\begin{equation*}
\inf_{\pi_T = \gamma du } I_{[0,T]} (\pi | \bar\rho)
\end{equation*}
measures the cost of joining $\bar\rho$ to $\gamma$ in the time
interval $[0,T]$ and 
\begin{equation}
\label{f06}
V(\gamma) \;:=\; \inf_{T>0} \inf_{\pi_T = \gamma du} 
I_{[0,T]} (\pi | \bar\rho)
\end{equation}
measures the cost of observing $\gamma$ starting from the stationary profile
$\bar\rho$. The functional $V$ is called the quasi-potential. 

It is expected in general that the quasi-potential equals the rate
functional of a large deviations principle for the empirical density
under the stationary state $\nu^N_{\alpha, \beta}$. We will see 
in Section~\ref{sec5} that this is indeed the case for the boundary
driven symmetric simple exclusion process. In particular, the
quasi-potential $V$ in \eqref{f06} is equal to the free energy $\mc F$ in
\eqref{f11}.

\section{Dynamical approach to stationary large deviations}
\label{sec5}

In this section we characterize the optimal path for the
variational problem \eqref{f06}. As a byproduct we show that $V$, as
defined in \eqref{f06}, equals $\mc F$, as defined in \eqref{f11}.
Unless explicitly stated, the arguments presented in this section hold
for interacting particle systems under general assumptions. We however
present - informally - the relevant statements and proofs in the context
of the one-dimensional boundary driven symmetric exclusion process. To
simplify the notation, given a density path $\pi\in D([0,T];\mc
M_+)$ such that $\pi_t$ is absolutely continuous with respect to the
Lebesgue measure for each $t\in [0,T]$, $\pi_t(du)= \lambda_t(u)du$,
we shall write $I_{[0,T]} (\lambda | \gamma)$ for $I_{[0,T]} (\pi |
\gamma)$.

\subsection{The reversible case.}
\label{sec5.1}

Let $\varphi_t(u)$ be the optimal path for the variational problem
\eqref{f06} on the interval $(-\infty,0]$ instead of $[0,\infty)$. 
In the reversible case, $\alpha = \beta$, from
Onsager-Machlup we expect that it is equal to the time reversal
of the relaxation trajectory $\rho_t(u)$ solution of \eqref{f04},
$\varphi_t(u) = \rho_{-t}(u)$. We show that this is indeed the case.

The cost of the path $\varphi$ is not difficult to compute. By
definition of $\varphi$ and by \eqref{f04}, $\partial_t \varphi_t =
-(1/2) \Delta \varphi_t$. In particular, $\nabla H_t = (\nabla
\varphi_t)/ \chi(\varphi_t)$ solves \eqref{f07} so that
\begin{equation*}
I_{(-\infty,0]} (\varphi | \bar\rho) \;=\;
\frac 12 \int_{-\infty}^0 dt \int_0^1 du\, \frac{(\nabla
\varphi_t)^2}{\chi(\varphi_t)}\;\cdot
\end{equation*}
Let $R(a) = \log \{a/1-a\}$. Rewrite the integrand as $\nabla
R(\varphi_t) \nabla \varphi_t = \nabla \{R(\varphi_t) - R(\bar\rho)\}
\nabla \varphi_t$, because $\bar\rho$ is constant in the reversible
case, and integrate by parts in space to obtain that
\begin{equation*}
I_{(-\infty,0]} (\varphi | \bar\rho) \;=\;
- \frac 12 \int_{-\infty}^0 \!dt \int_0^1 \!du\, \Delta \varphi_t 
\{R(\varphi_t) - R(\bar\rho)\} \;\cdot
\end{equation*}
Since $\partial_t \varphi = -(1/2) \Delta \varphi_t$, $\delta \mc F_0
(\varphi)/ \delta \varphi = R(\varphi) - R(\bar\rho)$, the previous
expression can be rewritten as
\begin{equation*}
\int_{-\infty}^0 \!dt \int_0^1 \!du\, \dot \varphi_t 
\frac{\delta \mc F_0 (\varphi_t)} {\delta \varphi_t}  
\;=\; \int_{-\infty}^0 \!dt  \, \frac d{dt} \mc F_0(\varphi_t) 
\;=\; \mc F_0(\varphi_0) - \mc F_0(\varphi_{-\infty}) 
\;=\; \mc F_0(\gamma)
\end{equation*}
because $\mc F_0(\bar\rho) =0$.  This proves that $V \le \mc F_0$.

The proof of Lemma \ref{s05} below, with $\nabla \{R(\lambda_t) -
R(\bar\rho)\}$ instead of $\nabla \{\delta W (\lambda_t) / \delta
\lambda_t\}$, shows that the cost of any trajectory $\lambda_t$
joining $\bar\rho$ to a profile $\gamma$ in the time interval $[0,T]$
is greater or equal to $\mc F_0(\gamma)$:
\begin{equation*}
I_{[0,T]} (\lambda | \bar\rho) \;\ge\; \mc F_0(\gamma)\;.
\end{equation*}
In particular, the trajectory $\varphi$ is optimal and $V(\gamma) =
\mc F_0(\gamma)$.

\subsection{The Hamilton-Jacobi equation}

We have seen in Subsection \ref{sec5.1} that the optimal path for
reversible systems is the relaxation path reversed in time.  In the
non reversible case, the problem is much more difficult and, in
general, we do not expect to find the solution in a closed form.  We
first derive a Hamilton-Jacobi equation for the quasi-potential by
interpreting the large deviation rate functional $I_{[0,T]} (\cdot |
\bar\rho)$ as an action functional
\begin{eqnarray*}
I_{[0,T]} (\lambda | \bar\rho) &=& \frac 12 \int_0^T dt \int_0^1 du\,
\frac 1{\chi(\lambda_t(u))} \Big\{ \nabla^{-1} \big[ \partial_t \lambda_t(u)
- (1/2) \Delta \lambda_t(u) \big] \Big\}^2 \\
&=:& \int_0^T dt\, \mc L(\dot \lambda_t, \lambda_t)\;.
\end{eqnarray*}
The quasi-potential $V$ may therefore be written as
\begin{equation}
\label{f08}
V(\gamma) \;=\; \inf_{T>0} \inf_{
\substack{ \lambda_0= \bar\rho \\
\lambda_T = \gamma }} \int_0^T\! dt\, \mc L(\dot \lambda_t,
\lambda_t)\;. 
\end{equation}
From this variational formula, taking the Legendre transform of the
Lagrangian, we derive the Hamilton-Jacobi equation for the
quasi-potential:
\begin{equation}
\label{HJ}
\Big\< \nabla \frac{\delta V(\gamma)}{\delta \gamma}, \chi(\gamma) 
\nabla \frac{\delta V(\gamma)}{\delta \gamma} \Big\> \; +\;
\Big\< \frac{\delta V(\gamma)}{\delta \gamma} , \Delta \gamma \Big\>
\;=\; 0
\end{equation}
and $\delta V(\gamma)/\delta \gamma$ vanishes at the boundary.

One is tempted to solve the Hamilton-Jacobi to find the
quasi-potential and then to look for a trajectory whose cost is given
by the quasi-potential. The problem is not that simple, however,
because the theory of infinite dimensional Hamilton-Jacobi equations
is not well established. Moreover, as well known, even in finite
dimension the solution may develop caustics in correspondence to the
Lagrangian singularities of the unstable manifold associated to the
stationary solution $\bar\rho$, see e.g.\ \cite{GT}.
Finally, the Hamilton-Jacobi equation has more than one solution.  In
particular, even if one is able to exhibit a solution, one still needs
to show that the candidate solves the variational problem \eqref{f08}.

The next lemma shows that a solution $W$ of the Hamilton-Jacobi
equation is always smaller or equal than the quasi-potential:

\begin{lemma}
\label{s05}
Let $W$ be a solution of the Hamilton-Jacobi equation \eqref{HJ}. Then,
$W(\gamma) - W(\bar\rho) \le V(\gamma)$ for all profiles $\gamma$.
\end{lemma}

\noindent{\bf Sketch of the proof.}
Fix $T>0$, a profile $\gamma$, and consider a path $\lambda$ in
$C_{\bar\rho}$ such that $\lambda_T= \gamma$. We need to show
that $I_{[0,T]} (\lambda | \bar\rho) \ge W(\gamma) - W(\bar\rho)$.

The functional $I_{[0,T]} (\lambda | \bar\rho)$ can be rewritten as
\begin{eqnarray}
\label{f09}
\!\!\!\!\!\!\!\!\!\!\!\!\!\!\! &&
\frac 12 \int_0^T dt \, \Bl \chi(\lambda_t) \, \Big\{ \nabla H_t - \nabla
\frac{\delta W (\lambda_t)}{\delta \lambda_t} \Big\}^2\Br \;+\;
\int_0^T dt \, \Bl \chi(\lambda_t) \, ( \nabla H_t) \, 
\Big( \nabla
\frac{\delta W (\lambda_t)}{\delta \lambda_t} \Big) \Br \nonumber \\
\!\!\!\!\!\!\!\!\!\!\!\!\!\!\! && \qquad
-\; \frac 12 \int_0^T dt \, \Bl \chi(\lambda_t) \, \Big\{ \nabla
\frac{\delta W (\lambda_t)}{\delta \lambda_t} \Big\}^2 \Br\;.
\end{eqnarray}
Since $\delta W (\lambda_t)/\delta \lambda_t$ vanishes at the
boundary, an integration by parts gives that the second integral
is equal to 
\begin{equation*}
- \int_0^T dt \, \Bl \frac{\delta W (\lambda_t)}{\delta \lambda_t} ,
\nabla \big( \chi(\lambda_t) \nabla H_t \big) \Br\;.
\end{equation*}
Since $W$ is a solution of the Hamilton-Jacobi equation, the third
integral is equal to
\begin{equation*}
\int_0^T dt \, \Bl \frac{\delta W (\lambda_t)}{\delta \lambda_t} ,
(1/2) \Delta \lambda_t \Br\;.
\end{equation*}
Summing this two expressions and keeping in mind that $H_t$ solves
\eqref{f07}, we obtain that $I_{[0,T]} (\lambda | \bar\rho)$ is
greater than or equal to 
\begin{equation*}
\int_0^T dt \, \Bl \frac{\delta W (\lambda_t)}{\delta \lambda_t} ,
\dot \lambda_t \Br \; =\; W(\lambda_T) - W(\lambda_0) \;=\;
W(\gamma) - W(\bar\rho) \;.
\end{equation*}
This proves the lemma.
\qed\smallskip

To get an identity in the previous lemma, we need the first term in
\eqref{f09} to vanish. This corresponds to have $\nabla H_t = \nabla 
\delta V (\lambda_t) / \delta \lambda_t$, i.e.\ to find a path
$\lambda$ which is the solution of
\begin{equation*}
\partial_t \lambda_t = (1/2) \Delta \lambda_t - \nabla
\Big \{ \chi(\lambda_t) \nabla 
\frac{\delta V (\lambda_t)}{\delta \lambda_t}\Big\}\;.
\end{equation*}
Its time reversal $\psi_t = \lambda_{-t}$, $t\in [-T,0]$ solves 
\begin{equation}
\label{f10}
\left\{
\begin{array}{l}
{\displaystyle \partial_t \psi_t = - (1/2) \Delta \psi_t
+ \nabla \Big \{ \chi(\psi_t) \nabla 
\frac{\delta V (\psi_t)}{\delta \psi_t}\Big\}\; ,} \\
{\displaystyle \psi_{-T} = \gamma \;, } \\
{\displaystyle \psi_t (0) = \alpha\;, 
\quad  \psi_t (1) = \beta\;. }
\end{array}
\right.
\end{equation}

As we argue in the next subsection, equation \eqref{f10} corresponds
to the hydrodynamic limit of the empirical density under the time
reversed dynamics; this is the Markov process on $\Sigma_N$ whose
generator is the adjoint to $L_N$ in $L_2(\Sigma_N,
\nu^N_{\alpha,\beta})$.

The next lemma shows that a weakly lower semi-continuous solution $W$
of the Hamilton-Jacobi equation is an upper bound for the
quasi-potential $V$ if one can prove that the solution of \eqref{f10}
relax to the stationary profile $\bar\rho$.

\begin{lemma}
\label{s06}
Let $W$ be a solution of the Hamilton-Jacobi equation \eqref{HJ},
lower semi-continuous for the weak topology. 
Fix a profile $\gamma$. Let  $\psi_t$ be the solution 
of \eqref{f10} with $V$ replaced by $W$.
If $\psi_0$ converges $\bar\rho$ for $T\uparrow\infty$,  
then $V(\gamma) \le W(\gamma) - W(\bar\rho)$.
\end{lemma}

\noindent{\bf Sketch of the proof.}  
To prove the lemma, given $\varepsilon>0$, it is enough to find
$T_\varepsilon>0$ and a path $\varphi_t$ such that $\varphi_0 =
\bar\rho$, $\varphi_{T_\varepsilon} = \gamma$, $I_{[0,T_\varepsilon]}
(\varphi| \bar\rho) \le W(\gamma) - W(\bar\rho)+ \varepsilon$.
Fix $T>0$ and let $\psi_t$ be the solution of equation \eqref{f10} in
the time interval $[-T,-1]$ with initial condition $\psi_{-T} =
\gamma$. Consider then an appropriate interpolation between
$\psi_{-1}$ and $\bar\rho$ which we again denote $\psi_t$,
$t\in[-1,0]$.  Let $\varphi_t = \psi_{-t}$, which is defined in the
time interval $[0,T]$. By definition of $I_{[0,T]}$,
\begin{equation*}
I_{[0,T]} (\varphi|\bar{\rho}) \;=\; I_{[0,1]} (\varphi|\bar{\rho})
\;+\; I_{[1, T]} (\varphi|\psi_{-1})\;.
\end{equation*}
Since $\psi_{-1}$ converges to $\bar\rho$ as $T\uparrow\infty$, the
first term can be made as small as we want by taking $T$
large. The second one, by definition of $\psi_t$ and by the
computations performed in the proof of Lemma \ref{s05}, is equal to
$W(\gamma) - W(\psi_{-1})$.  Since $\psi_{-1}$ converges to
$\bar\rho$ and since $W$ is lower semi-continuous we have
$W(\bar\rho)\le \liminf_{T\to\infty} W(\psi_{-1})$. Hence
$\limsup_{T\to\infty} I_{[0,T]}(\varphi|\bar\rho) \le W(\gamma) -
W(\bar\rho)$.  
This proves the lemma.
\qed

\smallskip
Putting together the two previous lemmata, we get the following statement.

\begin{theorem}
\label{s07}
Let $W$ be a solution of the Hamilton-Jacobi equation, lower
semi-continuous for the weak topology. 
Suppose that the solution $\psi_t$ of \eqref{f10}, with $V$ replaced
by $W$,  is such that  $\psi_0$ converges to $\bar\rho$ as 
$T\uparrow\infty$ for any initial profile $\gamma$. 
Then  $V(\gamma) = W(\gamma) - W(\bar\rho)$. 
Moreover, $\varphi_t = \psi_{-t}$ is the optimal path for the variational
problem \eqref{f08} defined in the interval $(-\infty,0]$ instead of
$[0,\infty)$.
\end{theorem}

\subsection{Adjoint hydrodynamic equation}
\label{s:ahe}

We have just seen that equation \eqref{f10} plays an important role in
the derivation of the quasi-potential. We show in this subsection that
\eqref{f10} describes in fact the evolution of the density profile
under the adjoint dynamics.

Consider a diffusive interacting particle system $\eta^N_t$  
satisfying the following assumptions.
The limiting evolution of the empirical density is described by a
differential equation
\begin{equation*}
\partial_t \rho = \mc D (\rho) \;,
\end{equation*}
where $\mc D$ is a differential operator. In the symmetric simple
exclusion process $\mc D(\rho)=(1/2)\Delta \rho $.  
Denote by $\xi^N_t = \eta^N_{-t}$ the time-reversed process.  The
limiting evolution of its empirical density is also described by a
differential equation 
\begin{equation}
\label{f16}
\partial_t \rho = \mc D^* (\rho) \;,
\end{equation}
for some integro-differential operator $\mc D^*$.  
Moreover the empirical densities satisfy a dynamical
large deviations principle with rate functions
\begin{equation*}
\frac 12 \int_0^T \!dt\,  
\Big\langle
\frac 1{\chi(\lambda_t)} 
\Big[ \nabla^{-1} \big(\partial_t \lambda_t
- \mc D^{i}( \lambda_t) \big) \Big]^2 
\Big\rangle
\;,\ \ \ i=1,2
\end{equation*}
where $\mc D^1=\mc D$ and $\mc D^2=\mc D^*$ for the original and 
time-reversed processes, respectively. In \cite{bdgjl1, bdgjl2} it is
shown that
\begin{equation}
\label{f15}
\mc D (\rho) \;+\; \mc D^* (\rho) \;=\; \nabla \Big( \chi(\rho) \nabla
\frac{\delta V}{\delta \rho} \Big)\;.
\end{equation}
In this general context, equation \eqref{f10} takes the form
\begin{equation*}
\partial_t \rho \;=\; - \mc D (\rho) + \nabla \Big( \chi(\rho) \nabla
\frac{\delta V}{\delta \rho} \Big) \;.
\end{equation*} 
Therefore, under the above assumptions on the dynamics, the solution
of \eqref{f10} represents the hydrodynamic limit of the empirical
density under the adjoint dynamics. In particular, for non reversible
systems, the typical path which creates a fluctuation is the
time-reverse of the relaxation path of the macroscopic dynamics of
$\xi^N_t$.  This principle extends the Onsager-Machlup theory to
irreversible systems.

\subsection{Optimal trajectory for the simple exclusion process}  

While all the arguments presented above are general, in this
subsection we obtain a more explicit description of the optimal trajectory 
for the variational problem \eqref{f08} for the boundary driven simple
exclusion process. As a corollary we show that the quasi-potential is
given by the expression \eqref{f11}. 

Let us first show how, in this case, it is possible to obtain a
solution of the Hamilton-Jacobi equation \eqref{HJ}. 
We look for a solution of the form
\begin{equation}\label{guess}
\frac{\delta W}{\delta\rho (u)} = 
\log \frac{\rho(u)}{1-\rho(u)} - \phi(u;\rho)
\end{equation}
for some functional $\phi(u; \rho)$ to be determined, satisfying the
boundary conditions $\phi(0 )= \log [\alpha /(1-\alpha)]$, 
$\phi(1)= \log [\beta /(1-\beta)]$.

Inserting \eqref{guess} into \eqref{HJ}, we get, note that 
$\rho - e^\phi/{(1+e^\phi)}$ vanishes at the boundary,
\begin{eqnarray*}
 0 &=& - \left\langle 
\nabla \left( \log \frac{\rho}{1-\rho} - \phi \right), \rho (1-\rho)
\nabla \phi \right\rangle
\\
&=&- \left\langle \nabla\rho, \nabla\phi\right\rangle 
+ \left\langle \rho (1-\rho), (\nabla\phi)^2 \right\rangle
\\
&=&-  \left\langle 
\nabla \left( \rho - \frac{e^\phi}{1+e^\phi} \right), \nabla \phi
\right\rangle
-\left\langle \left( \rho-\frac{e^\phi}{1+e^\phi}  \right)
\left(\rho-\frac{1}{1+e^\phi} \right), (\nabla\phi)^2 
\right\rangle
\\
&=& \left\langle\left( \rho - \frac{e^\phi}{1+e^\phi} \right) , 
\left( \Delta \phi + \frac{ (\nabla\phi)^2}{1+e^\phi} 
- \rho (\nabla\phi)^2 \right)\right\rangle
\end{eqnarray*}
We obtain a solution of the Hamilton--Jacobi equation if we solve the
following ordinary differential equation which relates the functional
$\phi(u)=\phi(u; \rho)$ to $\rho$
\begin{equation}\label{dphi}
\left\{
\begin{array}{l}{\displaystyle 
\frac {\Delta \phi(u)}{[\nabla\phi(u)]^2} + 
\frac{1}{1+e^{\phi(u)}} = \rho (u)  \quad u\in (0,1) }
\\
\\
\phi(0) = \log [\alpha/(1-\alpha)], 
\qquad \phi(1)= \log [\beta /(1-\beta)]
\end{array}
\right.
\end{equation}
A computation shows that the derivative of the functional 
\begin{equation}\label{Ssep}
W(\rho) = \int_{0}^1\!du \left\{
\rho\log \rho + (1-\rho)\log(1-\rho) + (1-\rho) \phi 
- \log \left(1+e^{\phi}\right) + \log
\frac{\nabla\phi}{\beta - \alpha} 
\right\}
\end{equation}
is given by \eqref{guess} when $\phi(u;\rho)$ solves \eqref{dphi}.
We note that by the change of variable $\phi = \log [ F /
(1-F)]$ equation \eqref{dphi} becomes \eqref{Deq} and 
\eqref{Ssep} becomes \eqref{f11}.

We next show that $W$ satisfies the hypotheses of Theorem~\ref{s07}.
In this case \eqref{f10} (with $V$ replaced by $W$), after shifting
the time interval from $[-T,0]$ to $[0,T]$, is given by 
the equation non local in space 
\begin{equation}
\label{adjHE}
\left\{
\begin{array}{l} 
{\displaystyle
\partial_t \rho  = \frac 12 \Delta \rho   
- \nabla \left\{ \rho  [ 1- \rho]  
\nabla \phi (\rho )
\right\} \quad u \in (0,1)
}
\\
{\displaystyle \rho_t(0) = \alpha, \qquad \rho_t (1)=\beta
} 
\\
{\displaystyle \rho_0(u)   = \gamma(u) } 
\end{array}
\right.
\end{equation}
where $\phi ( \rho )$ is to be obtained from $\rho $ by
solving \eqref{dphi}.  Since 
$\phi (\bar\rho) = \log [ \bar\rho/(1-\bar\rho)]$, 
we see that $\bar\rho$ is also a stationary solution of \eqref{adjHE}.

Equation \eqref{adjHE} can be related to the heat equation as follows.
Let $\rho_t$ be the solution of \eqref{adjHE} and introduce
$F=F_t(u)$ as
\begin{equation}
\label{Fphi}
F_t(u)=  \frac { e^{\phi(u;\rho_t )}}{1+e^{\phi(u;\rho_t )}} 
\end{equation}
it is not too difficult, see \cite[Appendix~B]{bdgjl2}, to check that
$F_t(u)$ satisfies the heat equation  
\begin{equation}\label{adjHEF}
\left\{
\begin{array}{l} 
{\displaystyle
\partial_t F_t (u) =  \frac 12 \Delta F_t (u)  
 \quad u \in (0,1)
}
\\
{\displaystyle F_t(0) = \alpha ,
\qquad F_t(1) = \beta,
} 
\\
{\displaystyle F_0(u)   =  
\frac { e^{\phi(u;\gamma )}}{1+e^{\phi(u;\gamma )}}}
\end{array}
\right.
\end{equation}
Conversely, given $F_t(u)$ which solves \eqref{adjHEF}, by setting 
$$
\rho_t (u) =  F_t(u) 
+  F_t(u) [1- F_t(u)]  \frac{ \Delta F_t(u)}{ [\nabla F_t(u) ]^2}
$$
a computation shows that $\rho_t $ solves \eqref{adjHE}. 

We have thus shown how a solution of the (non local, non linear)
equation \eqref{adjHE} can be obtained from the linear heat equation by
performing the non local transformation \eqref{Fphi} on the initial
datum.  In particular, since the solution $F_t(u)$ of \eqref{adjHEF}
converges as $t\to\infty$ to $\bar\rho$, we see that the functional
$W(\rho)$ given in \eqref{Ssep} satisfies the hypotheses of
Theorem~\ref{s07}. Since $W(\bar\rho)=0$ we thus have 
that the quasi-potential of the boundary driven simple exclusion
process is given by the functional \eqref{Ssep} where $\phi$ solves
\eqref{adjHE}.

To conclude the dynamical proof of Theorem~\ref{s02}, we have to 
identify the rate function $\mc F$  of the invariant measure with the
quasi-potential $V$.

\subsection{Free energy and quasi-potential}  

Bodineau and Giacomin \cite{bg}, adapting to this infinite dimensional
setting the method introduced by Freidlin and Wentzell \cite{FW} in the context
of small perturbations of dynamical systems, proved the following theorem 
which identifies $V$ and $\mc F$.

\begin{theorem}
\label{s11}
Let $I_{[0,T]}$ be the rate function in Theorem~\ref{s04} and define
the quasi-potential as in \eqref{f06}.  Then the empirical density
under the stationary state satisfies a large deviation principle with
rate functional given by the quasi-potential.
\end{theorem}

An explicit description of the quasi-potential, as the one here
discussed for the boundary driven simple exclusion process, is not
always possible. There are other few one dimensional boundary driven
models for which a similar representation has been obtained. 
This class includes the weakly asymmetric exclusion processes, the 
zero range processes and the KMP model \cite{bdgjl1,bgl,de}.

\section{Asymptotic behavior of the empirical current}
\label{sec4}

We examine in this section the current fluctuations over a fixed
macroscopic time interval.  In particular we discuss the law of large
numbers and the dynamical large deviations principle for the empirical
current. We state these results in the context of the boundary driven
symmetric exclusion process but similar results hold for more general
dynamics and for periodic boundary conditions.

Consider the boundary driven symmetric simple exclusion process
defined in Section \ref{sec1}.  For $0\le x\le N-1$, denote by
$j_{x,x+1}$ the rate at which a particle jumps from $x$
to $x+1$ minus the rate at which a particle jumps from $x+1$ to
$x$. For $x=0$, this is the rate at which a particle is created minus
the rate at which a particle leaves the system. A similar
interpretation holds at the right boundary.
An elementary computation shows that 
\begin{equation*}
j_{x,x+1} \;=\; N^2
\left\{
\begin{array}{ll}
\alpha - \eta(1) & \text{for $x=0$, } \\
\eta(x) - \eta (x+1) & \text{for $1\le x\le N-2$\;,} \\
\eta(N-1) - \beta & \text{for $x=N-1$. }
\end{array}
\right.
\end{equation*}
In view of \eqref{f01bis}, under the invariant measure
$\nu^N_{\alpha,\beta}$, the  average of $j_{x,x+1}$ is 
$$
E_{\nu^N_{\alpha, \beta}}[j_{x,x+1}] = N (\alpha - \beta)
$$

Given a bond $\{x,x+1\}$, $0\le x\le N-1$, let $J^{x,x+1}_t$ (resp.
$J^{x+1, x}_t$) be the number of particles that have jumped from $x$ to
$x+1$ (resp. $x+1$ to $x$) in the time interval $[0,t]$.  Here we
adopt the convention that $J^{0,1}_t$ is the number of particles
created at $1$ and that $J^{0,1}_t$ represents the number of
particles that left the system from $1$. A similar
convention is adopted at the right boundary.  The difference
$W^{x,x+1}_t= J^{x,x+1}_t - J^{x+1,x}_t$ is the net number of
particles flown across the bond $\{x, x+1\}$ in the time interval
$[0,t]$. 
Let us consider the stationary process $\bb P_{\nu^N_{\alpha,\beta}}$,
i.e.\ the boundary driven symmetric simple exclusion process in which
the initial condition is distributed according to the invariant
measure $\nu^N_{\alpha,\beta}$.
A simple martingale computation shows that $W^{x,x+1}_t /(Nt)$ 
converges, as $t\to\infty$, to $(\alpha - \beta)$ in probability. 
Namely, for each $N\ge 1 $, $x=0,\dots,N -1$, and  $\delta>0$, 
we have
\begin{equation*}
\lim_{t \to\infty} \bb P_{\nu^N_{\alpha, \beta}} 
\Big[ \, \Big\vert \frac{ W^{x,x+1}_t}{Nt} -  (\alpha - \beta) \Big\vert
\;>\; \delta \Big] \;=\; 0\;.
\end{equation*}

Let $\mc M$ be the space of bounded signed measures on $[0,1]$ endowed
with the weak topology.  For $t\ge 0$, define the \emph{empirical
integrated current} $ W^N_t \in \mc M$ as the finite signed measure on
$[0,1]$ induced by the net flow of particles in the time interval
$[0,t]$:
\begin{equation*}
W^{N}_{t} \;=\; N^{-2} \sum_{x=0}^{N-1}
W^{x,x+1}_t \delta_{x/N}\;.
\end{equation*}
Notice the extra factor $N^{-1}$ in the normalizing constant which
corresponds to the diffusive rescaling of time. In particular, for a
function $F$ in $C ([0,1])$, the integral of $F$ with respect to
$W^N_t$, also denoted by $\<W^N_t, F\> $, is given by
\begin{equation}
\label{empcurr}
\< W^N_t , F\> \;=\; N^{-2} \sum_{x=0}^{N-1}
F (x/N) \, W^{x,x+1}_t \;.
\end{equation}
It is not difficult to prove the law of large numbers for the
empirical current starting from an initial configuration associated to
a density profile.

\begin{proposition}
\label{s08}
Fix a profile $\gamma$ and consider a sequence of configurations
$\eta^N$ such that $\pi^N(\eta^N)$ converges to $\gamma (u) du$, as
$N\uparrow\infty$.  Let $\rho$ be the solution of the heat equation
\eqref{f04}. Then, for each $T>0$, $\delta>0$ and $F$ in $C([0,1])$,
\begin{equation*}
\lim_{N\to\infty} \bb P_{\eta^N} \Big[\, \Big\vert \<W^N_T, F\>
\;+\; (1/2) \int_0^T dt\, \int_{0}^1 F(u) 
\nabla \rho_t(u) \, du \Big\vert > \delta \Big ] \;=\; 0\; .
\end{equation*}
\end{proposition}

This result states that the empirical current $W^N_t$ converges to
the time integral of $- (1/2) \nabla \rho_t(u)$, which is the
instantaneous current associated to the profile $\rho_t$. Thus, if we
denote by $w(\gamma) = -(1/2) \nabla\gamma$ the instantaneous current
of a profile $\gamma$, we have that
\begin{equation*}
\lim_{N\to\infty} W^N_t \;=\; \int_0^t ds \, w(\rho_s)
\end{equation*}
in probability.
Proposition \ref{s08} is easy to understand. The local conservation of
the number of particles is expressed by
$$
\eta_t(x) - \eta_0(x) = W^{x-1,x}_t - W^{x,x+1}_t \;.
$$
It gives the following continuity equation for the empirical
density and current.  Let $G: [0,1] \to \bb R$ be a smooth function
vanishing at the boundary and let $(\nabla_N G) (x/N) = N\{ G(x+1/N) -
G(x/N)\}$. Then,
$$
\langle \pi^N_t, G \rangle - \langle \pi^N_0, G \rangle =
\langle W^N_t , \nabla_N G \rangle \;.
$$
The previous identity shows that the empirical density at time
$t$ can be recovered from the initial state and the empirical current
at time $t$.
In contrast, the empirical density at time $t$ and at time $0$
determines the empirical current at time $t$ only up to a constant.
Letting $N\uparrow\infty$ in the previous identity,  since $\pi^N$ 
converges to the solution of
the heat equation \eqref{f04}, an integration by parts gives that
\begin{equation*}
\langle W_t ,  \nabla G \rangle \;=\;
\langle \rho_t, G \rangle - \langle \rho_0, G \rangle \;=\;
\frac 12 \int_0^t ds\, \< \Delta \rho_s ,  G\> \;=\;
- \frac 12 \int_0^t ds\, \< \nabla \rho_s , \nabla G\> \;.
\end{equation*}
where $W_t$ is the limit of $W^N_t$.

After proving this law of large numbers for the current, we examine
its large deviations properties.  To state a large deviations
principle for the current we need to introduce some notation. Fix
$T>0$ and recall that we denote by $w(\gamma) = -(1/2) \nabla \gamma$
the instantaneous current associated $\gamma$. For a density profile
$\gamma$ and a path $W$ in $D([0,T], \mc M)$, denote by 
$w_t = \dot{W}_t$ and let $\rho^{\gamma, W}_t$ the weak solution of
\begin{equation}
\label{f17}
\left\{
\begin{array}{l}
\partial_t \rho_t \;+\; \nabla w_t \;=\; 0\;, \\
\rho_0 (u)  = \gamma (u)\;, \\
\rho_t(0) = \alpha\;, \quad \rho_t(1) = \beta \;.
\end{array}
\right.
\end{equation}
We note that the trajectory $\rho^{\gamma, W}_t$ is the one followed
by the density profile if the initial condition is $\gamma$ and the
instantaneous current is $w$.  As for the empirical density, the rate
functional for the empirical current is given by a variational
expression. Referring to \cite{bdgjl7} for the precise definition, we
here note that for trajectories $W$ in $D([0,T], \mc M)$ the rate
functional is finite only if the associated density path
$\rho^{\gamma, W}_t du$ belongs to $C([0,T], \mc M_+)$; moreover when
$W$ is a smooth path we have
\begin{equation}
\label{f12}
\mc I_{[0,T]} (W|\gamma) \; =\; \frac 12 \int_0^T dt \, 
\Big\< \frac 1{\chi(\rho^{\gamma, W}_t)}
\, \big\{ \dot{W}_t -  w(\rho^{\gamma, W}_t) \big \} ^2 \Big\>\; .
\end{equation}

The following theorem is proven in \cite{bdgjl7} in the case of
periodic boundary condition. The proof can easily be modified to cover
the present setting of the boundary driven simple exclusion process.

\begin{theorem}
\label{s09}
Fix $T>0$ and a smooth profile $\gamma$ bounded away from $0$ and $1$.
Consider a sequence $\eta^N$ of configurations associated to $\gamma$
in the sense that $\pi^N(\eta^N)$ converges to $\gamma(u) du$ as
$N\uparrow\infty$. Fix $W$ in $D([0,T], \mc M)$ and an associated 
neighborhood $V_\epsilon(W)$ of radius $\epsilon$.  Then,
\begin{eqnarray*}
\!\!\!\!\!\!\!\!\!\!\!\!\!\!\!\!\!\! &&
\limsup_{\epsilon\to 0} \limsup_{N\to\infty} \frac 1N \log
\bb P_{\eta^N} \big\{ W^N \in V_\epsilon(W) \big\}
\;\le\; - \mc I_{[0,T]} (W | \gamma) \; , \\
\!\!\!\!\!\!\!\!\!\!\!\!\!\!\!\!\!\! && \qquad
\liminf_{\epsilon\to 0} \liminf_{N\to\infty} \frac 1N \log
\bb P_{\eta^N} \big\{ W^N \in V_\epsilon(W) \big\}
\;\ge\; - \mc I_{[0,T]} (W | \gamma) \; .
\end{eqnarray*}
\end{theorem}

Since the trajectory of the empirical density can be recovered from
the evolution of the current and the initial condition, the large
deviations principle for the empirical density stated in
Theorem~\ref{s04} can be obtained from the large deviations principle
for the current by the  contraction principle, see \cite{bdgjl7} for
the proof.

\section{Large deviations of the time averaged empirical current}

In this section we investigate the large deviations properties of the 
mean empirical current $W^N_T/T$ as we let \emph{first} $N\to\infty$
and \emph{then} $T\to\infty$. As before, unless stated explicitly, 
the analysis carried out in this section does not depend on the
details of the symmetric simple exclusion process so that it holds in
a general setting.

Since the density is bounded, for $T$ large the time averaged
empirical current must be constant with respect to the space variable
$u$. This holds in the present one-dimensional setting; in higher
dimensions the condition required would be the vanishing of the
divergence. Indeed, if this condition were not satisfied we would have
an unbounded (either positive or negative) accumulation of particles.
We next discuss the asymptotic probability that the time averaged
empirical current equals some fixed constant.

For a smooth profile $\gamma$ bounded away from $0$ and $1$, let
$\tilde \Phi : \mc M \to [0,+\infty]$ be the functional defined by
\begin{equation}
\label{IT}
\tilde \Phi (J) = 
\begin{cases}
{\displaystyle 
\inf_{T>0} \frac 1T \, \inf_{W \in \mc A_{T,q} } \, 
\mc I_{[0,T]}  (W | \gamma) 
}
& \textrm{ if  $J(du) = q \,du$ \ for some \ $q\in\bb R$}  
\\
+\infty &\textrm{ otherwise}
\end{cases}
\end{equation}
where $\mc A_{T, q}$ stands for the set of currents with time average
equal to $q$
\begin{equation*}
\mc A_{T,q} := \Big\{ W \in D\big([0,T];\mc M \big) \,:\, 
 \frac 1T \int_0^T \!dt \, \dot{W}_t(du)  =  q \, du  \Big\}\; .
\end{equation*}

It is not difficult to show that $\tilde \Phi$ is convex. In the
present context of the boundary driven simple exclusion process, it is
also easy to verify that the functional $\tilde \Phi$ does not depend
on on the initial condition $\gamma$.  We emphasize however that, in
the case of periodic boundary condition $\tilde \Phi$ depends on
$\gamma$ only through its total mass $\int\!du\, \gamma(u)$.  Indeed,
we may start by driving the empirical density from a profile $\gamma$
to a profile $\gamma'$ in the time interval $[0,1]$ paying a finite
price, note that in the periodic case $\gamma$ and $\gamma'$ must have
the same mass. As $T\uparrow\infty$, this initial cost vanishes and the
problem is reduced to the original one starting from the profile
$\gamma'$.  Let us finally introduce $\Phi$ as the lower
semi-continuous envelope of $\tilde \Phi$, i.e.\ the largest lower
semi-continuous function below $\tilde \Phi$.  The next theorem states
that, as we let first $N\uparrow \infty$ and then $T\uparrow \infty$
the time averaged empirical current $W^N_T/T$ satisfies a large
deviation principle with rate function $\Phi$. We refer to
\cite{bdgjl7} for the proof which is carried out by analyzing the
variational problem $\inf_{W\in \mc A_{T,q}} T^{-1} \, \mc I_{[0,T]}
(W | \gamma)$ as $T\uparrow \infty$ and showing that it converges, in
a suitable sense, to the variational problem defining $\Phi$.

\begin{theorem}
\label{s10}
Fix $T>0$ and a smooth profile $\gamma$ bounded away from $0$ and $1$.
Consider a sequence $\eta^N$ of configurations associated to $\gamma$
in the sense that $\pi^N(\eta^N)$ converges to $\gamma(u) du$ as
$N\uparrow\infty$. Fix $J \in \mc M$ and a neighborhood
$V_\epsilon(J)$ of radius $\epsilon$.  Then,
\begin{eqnarray*}
&&
\limsup_{\epsilon\to 0} \limsup_{T\to\infty} \limsup_{N\to\infty}
\frac 1{T N} \log \bb P_{\eta^N}
\Big[ \frac 1T {W}^N_T \in V_\epsilon(J) \Big] \;\le\;
- \Phi (J) \;, \\
&&\qquad
\liminf_{\epsilon\to 0} \liminf_{T\to\infty} \liminf_{N\to\infty}
\frac 1{T N} \log \bb P_{\eta^N}
\Big[ \frac 1T {W}^N_T \in V_\epsilon(J) \Big] \;\ge\;
-  \Phi (J) \;.
\end{eqnarray*}
\end{theorem}
\medskip

A result analogous to Theorem~\ref{s10} can be proven for other diffusive
interacting particle systems. 
Consider a system with a weak external field $E=E(u)$,
whose hydrodynamic equation, describing 
the evolution of the empirical density on the macroscopic scale, 
has the form 
\begin{equation}
\label{f18}
\partial_t \rho_t = \nabla \big( D(\rho_t) \nabla\rho_t \big) - 
\nabla\big( \chi(\rho_t) E \big)\;.
\end{equation}
where $D(\rho)$ is the diffusion coefficient and $\chi(\rho)$ is the
mobility. For the symmetric simple exclusion process $D=1/2$ and
$\chi(\rho)=\rho(1-\rho)$.  In the general case, the large deviations
functional $\mc I_{[0,T]} (\cdot | \gamma)$ has the same form
\eqref{f12} with $w(\gamma)= - D(\gamma) \nabla \gamma + \chi(\gamma)
E$ and $\rho^{\gamma, W}$ the solution of \eqref{f17}.  For systems
with periodic boundary conditions, the boundary conditions in
\eqref{f17} is modified accordingly.  In the remaining part of this
section we analyze the variational problem \eqref{IT} for different
systems and show that different scenarios are possible.

A possible strategy for minimizing $\mc I_{[0,T]}(w|\gamma)$ with the
constraint that $w\in \mc A_{T,q}$, i.e.\ that the time average of $w$
is fixed, consists in driving the empirical density to a density
profile $\gamma^*$, remaining there most the time and forcing the
associated current to be equal to $q$.  This is the strategy
originally proposed by Bodineau and Derrida \cite{bd1}.  
In view of \eqref{f12} the asymptotic cost, as $T\uparrow\infty$,  
of this strategy is 
\begin{equation*}
\frac 12 \Big\langle \big[ q +  D(\gamma^*) \nabla \gamma^* \big], 
\frac 1{\chi (\gamma^*)} \big[ q +  D(\gamma^*) \nabla \gamma^* \big]
\Big\rangle\;.
\end{equation*}
If we minimize this quantity over all profiles $\gamma^*$ we obtain a
functional $U$ which gives the cost of keeping a current $q$ at a
fixed density profile:
\begin{equation}
\label{f14}
U(q) := \inf_{\rho}
\frac 12 \Big\langle \big[ q + D(\rho) \nabla \rho \big], 
\frac 1{\chi (\rho)} \big[ q + D(\rho) \nabla \rho \big]
\Big\rangle\;.
\end{equation}
where the infimum is carried out over all smooth density profiles
$\rho=\rho(u)$ bounded away from $0$ and $1$ which satisfy the
boundary conditions $\rho(0)=\alpha$, $\rho(1)=\beta$.  
As observed above, for boundary driven systems all density profiles
are allowed while for periodic boundary condition only profiles with
the same total mass $m=\int_0^1\!du \, \rho(u)$ are allowed. In the
latter case, the functional $U$ depends on the total mass $m$ and is
denoted by $U_m$.

As proven in \cite{bdgjl5,bdgjl6,bdgjl7}, for the symmetric simple exclusion
process the strategy above is the optimal one, i.e.\ $\Phi =U$. 
It is in fact not difficult to show that in this case $U$ is lower 
semi-continuous, so that $\tilde \Phi =\Phi$. More generally we have
the following result.  

\begin{lemma}
\label{s12}
Let $E=0$.  
If $D(\rho) \chi''(\rho) \le D'(\rho) \chi'(\rho)$ for any $\rho$,
then $\Phi=U$.
\end{lemma}

Besides the symmetric simple exclusion process, the hypothesis of the
lemma is also satisfied for the zero range model, where $D(\rho) =
\Psi '(\rho)$ and $\chi(\rho) =\Psi(\rho)$ for some strictly
increasing function $\Psi : \bb R_+ \to \bb R_+$, and for the non
interacting Ginzburg--Landau model, where $\rho\in\bb R$, $D(\rho)$ is
an arbitrary strictly positive function and $\chi(\rho)$ is constant.

For systems with periodic boundary condition we have shown
\cite{bdgjl6,bdgjl7} that the profile which minimizes $U_m$ is the
constant profile if $1/\chi(\rho)$ is convex.  

\begin{lemma}
\label{s13}
Let $E=0$. If the function $\rho\mapsto 1/\chi(\rho)$ is convex, then 
\begin{equation*}
U_m(q) \;=\; \frac 12 \frac{q^2}{\chi(m)}
\end{equation*}
and the constant profile $\rho(u)=m$ is optimal for the variational
problem \eqref{f14}.
\end{lemma}

The assumption of this lemma is satisfied by the symmetric simple
exclusion process as well as by the KMP model
\cite{bgl,kmp}, where $D(\rho)=1$ and $\chi(\rho) = \rho^2$.

\medskip 
As first discussed in \cite{bdgjl5}, the above strategy is not always
the optimal one, i.e.\ there are systems for which $\Phi < U$.  In
\cite{bdgjl5,bdgjl6} we interpreted this strict inequality as a
dynamical phase transition. In such a case the minimizers for
\eqref{IT} become in fact time dependent and the invariance under time
shifts is broken.  We now illustrate how different behaviors of the
variational problem \eqref{IT} leads to different dynamical regimes.
We consider the system in the ensemble defined by conditioning on the
event $ (T)^{-1} W^N_{T}(du) = q\, du$, $q\in \bb R$, with $N$ and
$T$ large.  The parameter $q$ plays therefore the role of an intensive
thermodynamic variable and the convexity of $\Phi$ expresses a
stability property with respect to variations of $q$.

If $\Phi(q)=U(q)$ and the minimum for \eqref{f14} is attained for
$\rho=\hat\rho(q)$ we have a state analogous to a unique phase: by
observing the system at any fixed time $t=O(T)$ we see, with probability
converging to one as $N,T\to \infty$, the density 
$\pi^N_T \sim\hat\rho(q)$ and the instantaneous current 
$ {\dot W}^N_t \sim q$.

While the functional $\Phi$ is always convex, $U$ may be not convex;
an example of a system with this property is given in \cite{bdgjl6}.  If
$\Phi$ is equal to the convex envelope of $U$, we have a state
analogous to a phase coexistence. Suppose for example $q= p q_1 +
(1-p)q_2$ and $U(q) > U^{**} (q)= p U(q_1) + (1-p) U(q_2)$ for some
$p,q_1,q_2$; here $U^{**}$ denotes the convex envelope of $U$.
The values $p,q_1,q_2$ are determined by $q$ and $U$.
The density profile is then not determined, but rather we observe with
probability $p$ the profile $\hat\rho(q_1)$ and with probability $1-p$
the profile $\hat\rho(q_2)$. 

Consider now the case in which a minimizer for \eqref{IT} is a 
current path $w_t$ not constant in $t$ and  
suppose that it is periodic with period $\tau$. We denote by 
$\hat\rho_t$ the corresponding density. Of course we have
$\tau^{-1}\int_0^\tau\!dt \, w_t = q$. In this situation we
have in fact a one parameter family of minimizers which are obtained
by a time shift $\alpha\in[0,\tau]$.  
This behavior is analogous to a non translation invariant state in equilibrium
statistical mechanics, like a crystal. Finally, if the optimal path
for \eqref{IT} is time dependent and not periodic the corresponding
state is analogous to a quasi-crystal.

\medskip
The explicit formula for $U_m$ derived in Lemma~\ref{s13} permits to
show that under additional conditions on the transport coefficient $D$
and $\chi$, a dynamical phase transition occurs. We discuss only the
case of periodic boundary conditions.
In this situation a time-averaged current $q$ may be produced 
using a traveling wave density profile, $\rho_t(u) =\rho_0(u-vt)$,
with velocity $v \sim q$. 
Assume now that $E=0$ and the function $\rho\mapsto \chi(\rho)$ 
is strictly convex for $\rho=m$. 
Then, for sufficiently large $q$, the traveling wave strategy is more
convenient than the one using the constant profile $m$ \cite{bdgjl5,bdgjl6}.
In particular, if $\rho\mapsto 1/\chi(\rho)$ is convex so that Lemma~\ref{s13}
can be applied, we have 
\begin{equation}
\label{ds}
\Phi_m(q) \;<\; U_m(q)
\end{equation}
for sufficiently large $q$.  In the KMP model the above hypotheses are
satisfied for any $m>0$; we can thus conclude that a dynamical phase
transition takes place for sufficiently large time-averaged currents.

The above analysis can also be applied to the weakly asymmetric simple
exclusion process \cite{bdgjl6}. It yields that if $|E/q|>
[m(1-m)]^{-1}$ for $q$ large there exists a traveling wave whose cost
is strictly less than the one of the constant profile $\rho(u) = m$.
The analysis in \cite{bd2} suggest however that the strict inequality
\eqref{ds} holds also in this case.
Moreover, the numerical simulations in \cite{bd2} indicate the
existence, for the weakly asymmetric simple exclusion process, of a
critical current $q^*$ below which the optimal profile is constant and
above which the optimal profile is a traveling wave.

\bigskip
\noindent{\bf Acknowledgments.}
The authors acknowledge the support of PRIN MIUR 2004\-028108 and
2004015228.  A.D.S.\ was partially supported by Istituto Nazionale di
Alta Matematica.  C.L.\ acknowledges the partial support of the John
S.\ Guggenheim Memorial Foundation, FAPERJ and CNPq.

\end{document}